\documentclass[12pt]{article}%
\usepackage{amssymb}
\usepackage{amsfonts}
\usepackage{amsmath}
\usepackage{graphicx}%
\begin{document}

\begin{center}
{\Large Huffman coding as an algorithm to construct chains }

{\Large in partition lattices}

\bigskip

Stephan Foldes

Tampere University of Technology

33101 Tampere,\ Finland

sf@tut.fi

\bigskip

\bigskip\textbf{Abstract}
\end{center}

\textit{The} \textit{Huffman coding algorithm is interpreted in the context of
the lattice of partitions of the source alphabet. Maximal chains in the
partition lattice correspond to linear extensions of tree orders, and those
among the chains that exhibit a simple greedy property correspond precisely to
executions of the Huffman algorithm.}

\bigskip

\bigskip

The input to the Huffman algorithm of information theory is a non-empty finite
set $S,$ called \textit{source} \textit{alphabet, }together with a map\textit{
}$p$ associating to each \textit{source symbol} $i\in S$ a non-negative real
number $p(i)$. For each subset $A\subseteq S$ we also write $p(A)$ for $%
{\displaystyle\sum\limits_{i\in A}}
p(i).$

The Huffman algorithm constructs a binary tree whose terminal nodes are the
source symbols, source symbol $i$ being at distance $l_{i}$ in the tree, such
that the sum%
\[
\sum_{i\in S}p(i)l_{i}\text{
\ \ \ \ \ \ \ \ \ \ \ \ \ \ \ \ \ \ \ \ \ \ \ \ \ \ \ \ \ \ \ \ \ \ \ \ \ (1)}%
\]
is minimal among all binary trees with terminal node set $S.$ (A binary
Huffman code is then produced by topologizing the tree, designating the sons
of each non-terminal node as left or right son, and recording the left-right
direction sequence of the path to each terminal node from the root as a 0-1
sequence.) The Huffman algorithm is non-deterministic, in the sense that
arbitrary choices made during its execution may produce different, even
non-isomorphic, but equally optimal trees.

It is useful to consider also the sum (1) not only for binary trees, but also
for binary forests with terminal node set $S$ (where $l_{i}$ is the distance
of $i$ from the root of the tree component to which $i$ belongs). Note that we
omit the usual assumption that the sum of the various $p(i)$ be $1$, as this
plays a role only in the interpretation of the input and the output of the
Huffman algorithm, but not in the procedure itself. The Huffman algorithm
constructs the optimal binary tree by building a sequence of binary forests
with the same terminal node set $S$, starting with the forest in which every
tree component has only one node (a source symbol), moving at each step from a
given forest to a new one with one less tree component, by including two of
the tree components $T_{1},T_{2}$ in a single new tree component $T$\ so that
the root of $T$ (the only new node in the next forest) will have the roots of
$T_{1}$ and $T_{2\text{ }}$as its sons. Clearly the end result is a tree, and
the number of forests involved in the sequence is the number of source
symbols. The specificity of the Hufmann algorithm, which turns out to
guarantee the optimality of the tree conctructed, is that at each step the
tree components $T_{1}$ and $T_{2\text{ }}$are chosen so that the sum of the
numbers $p(i)$ taken over the terminal nodes of the two chosen tree components
is minimal, among all possible choices of two tree components.

It is irrelevant what the nodes of the forests actually are, therefore to have
a canonical representation of the forests we shall suppose that each terminal
node is a singleton containg a source symbol, and each node \textit{is }simply
the set of those source symbols at the terminal nodes that are its descendants
(which set does not change from the stage when the node is added until the
procedure terminates). Then the node set of the tree resulting at the end of
an execution of the Huffman algorithm is nothing else, in the terminology of
comparable-or-disjoint set families (see $\left[  \text{CS},\text{ CHS,
HR}\right]  $), than a CD basis minus the empty set, in the lattice of subsets
of the source alphabet. At each stage of the algorithm what we have is a
CD-independent set $F$\ of non-empty sets of source symbols, such that for
each maximal member $X$ of $F,$ the set
\[
\left\{  Y\in F:Y\subseteq X\right\}
\]
together with $\emptyset$\ is a CD bases of the lattice of subsets of $X.$
This is in fact a complete characterization of the Huffman algorithm. The
algorithm can also be described as follows:

\bigskip

\textbf{Huffman Algorithm} \ \textit{Given a non-empty finite set} $S$
\textit{of} $n$ \textit{elements}, \textit{with a non-negative real number}
$p(i)$ \textit{assigned to each} $i\in S.$ \textit{Let} $\mathcal{C}$
\textit{be the set of maximal chains in the lattice of partitions of a finite
non-empty set} $S$, \textit{construct a sequence} $F_{1}\subset F_{2}%
\subset...$ \textit{of families of subsets of} $S$, \textit{such that}

(i) $F_{0}$ \textit{is the family of all singletons, }

(ii) \textit{while} $S\notin F_{k}$, \textit{let} $F_{k+1}$ \textit{consist of
the members of} $F_{k}$ \textit{and of the union of two distinct maximal
members} $A\cup B$ \textit{of} $F_{k}$ \textit{such that} \textit{the sum}%

\[
\sum_{i\in A\cup B}p(i)l_{i}%
\]
\textit{is of maximal value over all choices of distinct maximal members}
$A,B$ \textit{of} $F_{k}$ .

\textit{Then the sequence stops at stage} $n-1$, $F_{0}\subset F_{1}%
\subset...\subset F_{n-1}$ \textit{where} $F_{n-1}\cup\left\{  \emptyset
\right\}  $ \textit{is a CD basis in the lattice of subsets of} $S,$
\textit{minimizing over all CD bases} $F$ \textit{the sum}%
\[
\sum_{i\in S}p(i)\text{ }\left[  Card\text{ }\left\{  A\in F:A\neq S\text{,
}i\in A\right\}  \right]  =%
{\displaystyle\sum\limits_{\substack{A\in F\\A\neq S}}}
p(A)
\]

\bigskip

Given any finite non-empty set $S,$ say with $n$ elements, consider a CD basis
$F$ of the lattice of subsets of $S.$\ Then the non-empty members of $F$
constitute a binary tree $T$ with root $S$, the terminal nodes of which are
the singleton subsets of $S$. Consider the set $F^{\prime}$ of non-singleton
members of $F$, and any linear extension\ $\lambda$\ of the subset relation
(partial order) of $F^{\prime}.$ Clearly $F^{\prime}$ has $n-1$ members, let
these be enumerated $A_{1},...,A_{n-1}$ in the increasing order of $\lambda.$
Then to $\lambda$ we can associate a maximal chain $\Pi_{0}<\Pi_{1}%
<...<\Pi_{n-1}$ in the lattice of partitions of $S$, where for $1\leq k\leq
n-1$, the partition $\Pi_{k}$ is obtained from $\Pi_{k-1}$ by replacing with
$A_{k}$ those two (uniquely determined) partition classes of $\Pi_{k-1}$ the
union of which is $A_{k}.$ In fact all maximal chains in the partition lattice
can be obtained this way, which we can also state in terms of the reverse
construction, in Proposition 1 below. Let $\mathcal{C}$ denote the set of
maximal chains in the lattice of partitions of a finite non-empty set $S.$ For
each maximal chain $K$ of partitions $\Pi_{0}<\Pi_{1}<...$ of $S$ and each
non-singleton set $A$ appearing as a class in any member of $K$, call the
positive integer $\min\left\{  k:0<k,\text{ }A\text{ is a class of }\Pi
_{k}\right\}  $ the \textit{index of first appearance} of $A$ in $K$. Note
that no two different sets can have the same index of first appearance, and
thus the non-singleton classes of the various members of $K$ are linearly
ordered according to their index of first appearance.

\bigskip

\textbf{Proposition} \ \textit{Let} $\mathcal{C}$ \textit{be the set of
maximal chains in the lattice of partitions of a finite non-empty set} $S.$
\textit{For each CD basis} $F$ \textit{of the lattice of subsets of} $S$,
\textit{let} $\Lambda_{F}$ \textit{be the set of linear extensions of the set
containment order on the set} $F^{\prime}$ \textit{of\ non-empty,
non-singleton members of} $F,$ and let $\Lambda$ \textit{be the union of all
the} $\Lambda_{F}$ \textit{for the various CD bases} $F$.

\textit{Then the map, which to each maximal chain} $K$ \textit{of partitions}
$\Pi_{0}<\Pi_{1}<...$ \textit{of} $S$ \textit{associates the set of
non-singleton classes of the various members of} $K,$\ \textit{linearly
ordered according to their indices of first appearance in} $K,$i\textit{s a
bijection from} $\mathcal{C}$ \textit{to} $\Lambda.$

\bigskip

\textbf{Greedy Chain Characterisation} \ \textit{The sequence }$F_{0}\subset
F_{1}\subset...\subset F_{n-1}$\textit{\ of sets of source symbols constructed
by any execution of the Huffman Algorithm on an }$n$\textit{-element finite,
non-empty set }$S,$ \textit{with a non-negative real number} $p(i)$
\textit{assigned to each} $i\in S,$ \textit{defines a maximal chain} $\Pi
_{0}<\Pi_{1}<...<\Pi_{n-1}$ \textit{in the lattice of partitions of }$S$,
\textit{where for} $1\leq k\leq n-1$, \textit{the partition} $\Pi_{k}$
\textit{is obtained from} $\Pi_{k-1}$ \textit{by replacing with the only set}
$C\in$ $F_{k}\setminus F_{k-1}$ \textit{the two classes} $A,B$ \textit{of}
$\Pi_{k-1}$ \textit{for which} $C=A\cup B.$

\textit{A maximal chain }$\Pi_{0}<\Pi_{1}<...<\Pi_{n-1}$\textit{ in the
lattice of partitions of a finite non-empty set} $S$ \textit{corresponds to an
execution of the Huffman Algorithm if and only if for each} $1\leq k\leq n-1$
\textit{the (uniquely determined) clas}s $C$ $\in$ $\Pi_{k}\setminus\Pi_{k-1}$
\textit{satisfies} $p(C)\leq p(K),$ \textit{for the sets} $K\in\Pi$
$\setminus$ $\Pi_{k-1}$ \textit{for all partitions} $\Pi$ \textit{covering}
$\Pi_{k-1}$ \textit{in the partition lattice}.

\bigskip

Acknowledgements.

This work has been co-funded by Marie Curie Actions and supported by the
National Development Agency (NDA) of Hungary and the Hungarian Scientific
Research Fund (OTKA, contract number 84593), within a project hosted by the
University of Miskolc, Department of Analysis.

The author wishes to thank S\'{a}ndor Radeleczki for useful comments and discussions.

\bigskip

\includegraphics[height=15mm, width=20mm]{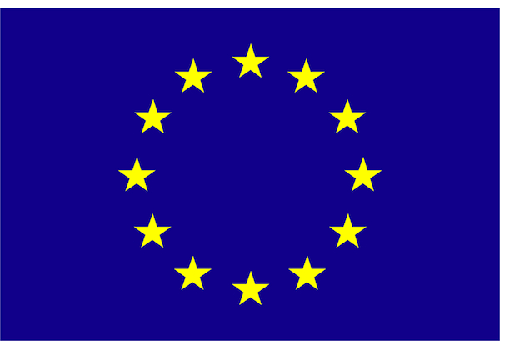}
\includegraphics[height=15mm, width=20mm]{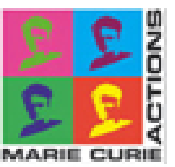}
\includegraphics[height=15mm, width=20mm]{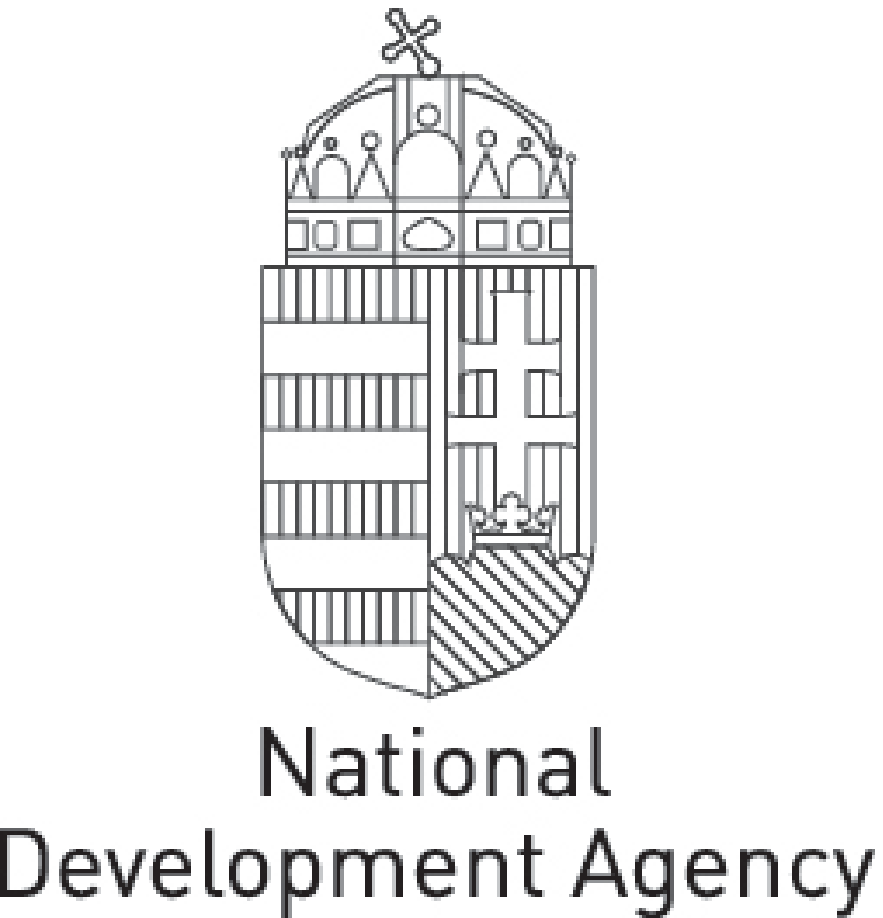} \includegraphics[height=15mm, width=20mm]{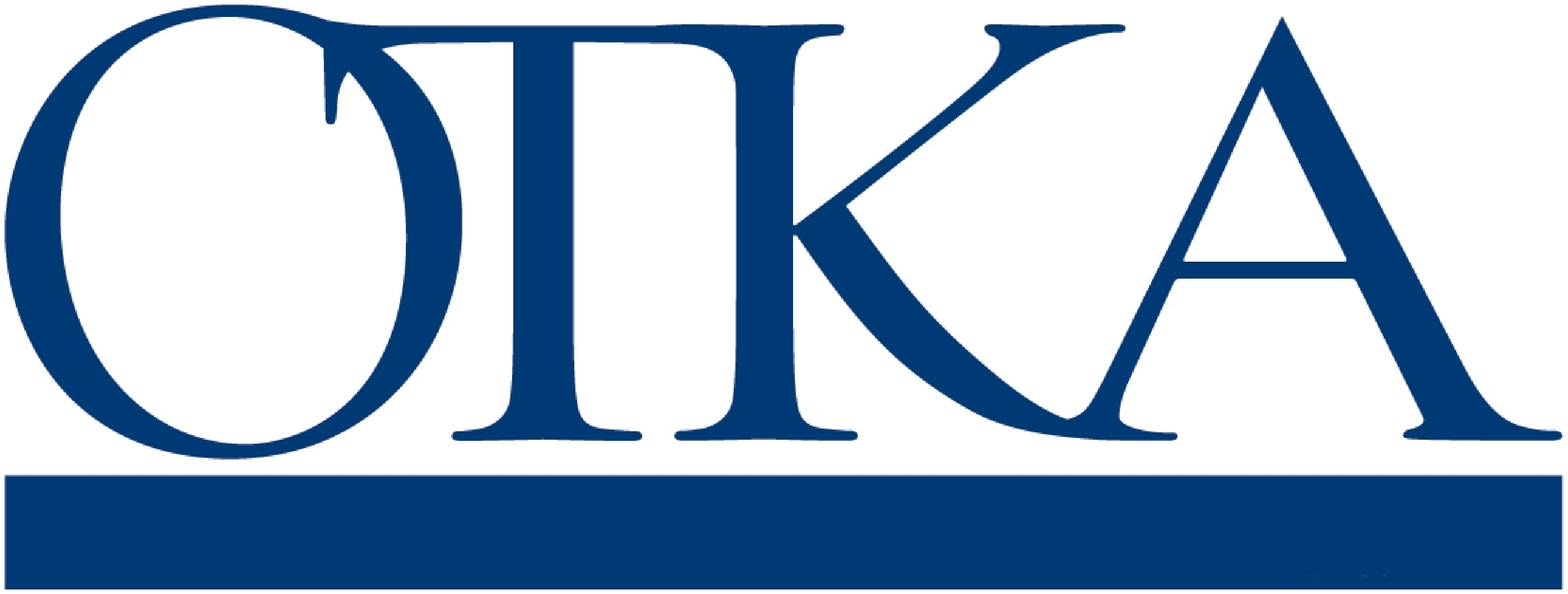}

\bigskip

\textbf{References}

\bigskip

$\left[  \text{CS}\right]  $ G. Cz\'{e}dli, E. T. Schmidt, CDW-independent
subsets in distributive lattices, Acta Sci. Math. (Szeged) 75 (2009), 49-53.

\bigskip

$\left[  \text{CHS}\right]  $ G. Cz\'{e}dli, M. Hartmann and E.T. Schmidt:
CD-independent subsets in distributive lattices, Publicationes Mathematicae
Debrecen, 74/1-2 (2009), 127-134

\bigskip

$\left[  \text{HR}\right]  $ E.K. Horv\'{a}th, S. Radeleczki, Notes on
CD-independent subsets, Acta Sci. Math. (Szeged) 78 (2012), 3-24

\end{document}